\documentclass[a4paper,12pt]{article}
\usepackage[top=2.5cm,bottom=2.5cm,left=2.5cm,right=2.5cm]{geometry}
\usepackage{amsfonts}
\usepackage{latexsym,amsmath,amssymb,amsfonts,epsfig,psfrag,url,graphics,ifpdf,multicol}
\usepackage{color,colordvi,amscd,amsthm} %eepic
\usepackage{tikz}
\usepackage[numbers,sort&compress]{natbib}
\usepackage{indentfirst,graphics,epsfig}
\usepackage{graphicx}
\usepackage{graphics}
\usepackage{graphicx,psfrag}
\usepackage{graphics}
\usepackage{mathrsfs}
\usepackage{tabularx}
\usepackage{booktabs}
\usepackage{floatrow}
\usepackage{multirow}
\usepackage{graphicx}
\usepackage{caption}
\usepackage[misc]{ifsym}
\usepackage[ruled]{algorithm2e}
\newfloatcommand{capbtabbox}{table}[][\FBwidth]

\newtheorem{theo}{Theorem}[section]
\newtheorem{lem}[theo]{Lemma}
\newtheorem{obs}[theo]{Observation}

\newtheorem{coro}[theo]{Corollary}

 \theoremstyle{definition}

\theoremstyle{remark}

\def\pf{\noindent{\bf Proof.\ }}
\allowdisplaybreaks

\usepackage{indentfirst,subfig}
\begin{document}
 %\thispagestyle{empty}
    %\rule{0cm}{0.5mm}
    \captionsetup[figure]{labelfont={bf},name={Fig.},labelsep=period}

\begin{center} {\large
Polynomial time recognition of vertices contained in all (or no) maximum dissociation sets of a tree}
\end{center}
\pagestyle{empty}

\begin{center}
	{
		{\small Jianhua Tu$^1$, Lei Zhang$^2$, Junfeng Du$^{3,}$, Rongling Lang$^{4,}$\footnote{Corresponding author.\\\indent \ \  E-mail: tujh81@163.com (J. Tu); 2018200896@mail.buct.edu.cn (L. Zhang); dujf1990@163.com (J. Du); ronglinglang@163.com (R. L)}}\\[2mm]
			
		{\small $^1$ School of Mathematics and Statistics, Beijing Technology and Business University, \\
			\hspace*{1pt} Beijing, P.R. China 100048} \\
		{\small $^2$ School of Mathematics and Statistics, Beijing Institute of Technology, \\
			\hspace*{1pt} Beijing, P.R. China 100081} \\	
		{\small $^3$ Department of mathematics, Beijing University of Chemical Technology, \\
			\hspace*{1pt} Beijing 100029, China} \\
		{\small $^4$ School of Electronics and information Engineering, Beihang University, \\
			\hspace*{1pt} Beijing, P.R. China 100191} \\
	}
	
\end{center}

\begin{center}
\begin{abstract}
In a graph $G$, a dissociation set is a subset of vertices which induces a subgraph with vertex degree at most 1. Finding a dissociation set of maximum cardinality in a graph is NP-hard even for bipartite graphs and is called the maximum dissociation set problem. The complexity of maximum dissociation set problem in various sub-classes of graphs has been extensively studied in the literature. 
In this paper, we study the maximum dissociation problem from different perspectives and characterize the vertices belonging to all maximum dissociation sets, and to no maximum dissociation set of a tree. We present a linear time recognition algorithm which can determine
whether a given vertex in a tree is contained in all (or no) maximum
dissociation sets of the tree. Thus for a tree with $n$ vertices, we can find all vertices
belonging to all (or no) maximum dissociation sets of the tree in $O(n^2)$ time.\vspace{5mm}



\end{abstract}
%\end{minipage}
\end{center}

\baselineskip=0.24in
\section{Introduction}

We consider only simple and undirected labeled graphs, and follow the terminology and notation of 
\cite{Bondy2008}. Let $G=(V,E)$ be a graph and $v$ be a vertex of $G$, we write $N_G(v)$ to denote the (open) neighborhood of $v$. 
The closed neighborhood of $v$ is defined as $N_G[v]=N_G(v)\cup\{v\}$.
The degree of $v$ is defined as $d_G(v)=|N_{G}(v)|$. For $U\subseteq V$, we write $G[U]$ to denote the subgraph induced by $U$. The subgraph $G[V\setminus U]$ is denoted by $G-U$. Furthermore, $G-U$ can be written by $G- u$ if $U = \{u\}$.

In a graph $G$, an independent set is a subset of vertices spanning no edges. Finding an independent set of maximum cardinality in a graph is a widely studied well-known problem of graph theory and is called the maximum independent set problem. In 1982, Hammer et al. \cite{Hammer1982} studied the maximum independent problem from different perspectives and investigated the vertices contained in all or in no maximum independent sets of a graph. Since then, researchers have extensively stuied this kind of problem for some other vertex subsets with
given properties. For example, Mynhardt \cite{Mynhardt1999}, Cockayne et al. \cite{Cockayne2003},
and Blidia et al. \cite{Blidia2005} considered this kind of problem for minimum dominating sets, total dominating sets, and
minimum double dominating sets of trees, respectively. Recently, Bouquet et al. \cite{Bouquet2021} studied this kinde of problem of minimum dominating sets on claw-free graphs, chordal graphs and triangle-free graphs.

In a graph $G$, a \emph{dissociation set} is a subset of vertices $F$
such that the induced subgraph $G[F]$ has maximum degree at most 1. A
\emph{maximum dissociation set} of $G$ is a dissociation set of
maximum cardinality. The \emph{dissociation number} $\psi(G)$ of a graph $G$ is the cardinality of a maximum dissociation
set of $G$. The concept of dissociation set was introduced by Yannakakis
\cite{Yannakakis1981} in 1981 and is a natrual generalization of independent set. The problem of finding a maximum dissociation set in a graph has been extensively studied on various sub-classes of graphs \cite{Alekseev2007,Bresar2011,Cameron2006,Orlovich2011} and is NP-hard for bipartite graphs \cite{Yannakakis1981}. 

A $k$-path vertex cover in a graph $G$ is a subset of vertices intersecting every $k$-path of $G$, where a $k$-path is a path of order $k$. It is easy to see that a set $S$ of vertices of a graph $G$ is a 3-path vertex cover of $G$ if and only if its complement $V(G)\setminus S$ is a dissociation set of $G$.
In this decade, the problem of finding a minimum
$k$-path vertex cover in a graph has received great attention \cite{Bresar2011,Kardos2011,Katrenic2016,Tu2011,Xiao2017}.

The main purpose of this paper is to characterize the vertices contained in all maximum dissociation sets and in no maximum dissociation set of a tree. Define the vertex subsets $\mathcal{A}(G)$,
$\mathcal{F}(G)$ and $\mathcal{N}(G)$ by
\begin{align}
\mathcal{A}(G)&=\{v\in  V(G):\text{$v$ is in all maximum dissociation sets of $G$}\}, \notag\\
\mathcal{F}(G)&=\{v\in  V(G):\text{$v$ is in some but not all
maximum dissociation sets of $G$}\}, \notag\\
\mathcal{N}(G)&=\{v\in  V(G):\text{$v$ is in no maximum dissociation
set of $G$}\}\notag.
\end{align}\par

On the other hand, the study is also
inspired by the relationship between the characteristic of vertex
subsets $\mathcal{A}(G)$ and $\mathcal{N}(G)$ and the number of
maximum dissociation sets in a graph $G$. In \cite{Tu2020}, Tu,
Zhang and Shi found four structure theorems concerning the vertex
subsets $\mathcal{A}(T)$ and $\mathcal{N}(T)$ for a tree $T$ and determined the
maximum number of maximum dissociation sets in a tree of order $n$.

The paper is organized as follows. In Section 2, we introduce some necessary
notation and lemmas. In Section 3, we characterize the vertex
subsets $\mathcal{A}(T)$ and $\mathcal{N}(T)$ of a tree $T$. In Section 4, a linear recognition algorithm which
can determine whether a given vertex in a tree is contained in all (or no)
maximum dissociation sets of the tree will be presented. Thus, using
the recognition algorithm, all vertices contained in all (or no)
maximum dissociation sets of a tree of order $n$ can be found in
$O(n^2)$ time.

\section{Preliminary results}

A \textit{rooted} tree $T$ is a connected acyclic graph with a specified vertex $r$,
called the \textit{root} of $T$. 
%An orientation of a rooted tree in which every vertex but the root has indegree one. 
Let $T$ be a tree rooted at $r$ and $v$ be a vertex of $T$, each vertex on the path from the root $r$ to the vertex $v$, including the
vertex $v$ itself, is called an \textit{ancestor} of $v$, and a \textit{descendant} of $v$ is a vertex $u$ such that $v$ is an ancestor of $u$.
An \textit{ancestor}  or \textit{descendant}  of a vertex is proper
if it is not the vertex itself. The \textit{parent} $p(v)$ of $v$ is the immediate proper ancestor of $v$, a \textit{child} of $v$ is a vertex $u$ such that $p(u)=v$. Define the vertex subsets $C_T(v)$, $D_T(v)$ and
$D_T[v]$ by
\begin{equation}
C_T(v)=\{u\in V(T) :\text{$u$ is a child of $v$}\}, \notag
\end{equation}
%\begin{equation}
%c(v)=u \quad \text{if}\quad C(v)=\{u\},\notag
%\end{equation}
\begin{equation}
D_T(v)=\{u\in V(T) :\text{$u$ is a proper descendant of
$v$}\},\notag
\end{equation}
\begin{equation}
D_T[v]=D(v)\cup\{v\}. \notag
\end{equation}
If no confusion occurs, these also be written by $C(v)$, $D(v)$ and
$D[v]$, respectively. We write $T_v$ to denote the subtree induced by $D_T[v]$.\par

A leaf in a tree is a vertex with degree 1, a
branch vertex is a vertex with degree at least 3. We write $B(T)$ to denote the set of branch
vertices of $T$. A path $P$ in $T$ is called a $v-L$ path, if $P$ joins $v$ to a leaf of $T$. Denote the order
of $P$ by $n(P)$, and for $i=0,1,2$, define
\begin{equation}
C^i(v)=\{u\in C(v):\text{$T_u$ contains a $u-L$ path $P$ with
$n(P)\equiv i\bmod3$}\} \notag
\end{equation}\par
Now, some basic observations about maximum dissociation sets
of the path are given.

\begin{obs}\label{obs2.1}
Let $P_n$ be a path of order $n\ge3$ and $u,v$ be the
leaves of $P_n$.\par 
(a) $\psi(P_n)=\tfrac{2n+i}{3}$, where $n\equiv
i(\bmod3)$, $i=0,1,2$.\par 
(b) If $n\equiv 0(\bmod3)$, then there exists a maximum dissociation set of $P_n$ that contains exactly one leaf.\par
(c) If $n\equiv 1(\bmod3)$, then both leaves of $P_n$ belong to all maximum dissociation sets of $P_n$, furthermore, $P_n$ has a
maximum dissociation set $F$ such that $d_{P_n[F]}(u)=0$.
\par 
(d) If $n\equiv 2(\bmod 3)$, then there is only one maximum
dissociation set $F$ in $P_n$, furthermore, $\{u,v\}\subset F$ and
$d_{P_n[F]}(u)=d_{P_n[F]}(v)=1$.
\end{obs}

Firstly, we characterize $\mathcal{A}(T)$ and $\mathcal{N}(T)$ in the case where $B(T)\le1$.

\begin{lem}\label{lem2.1}
Let $T$ be a rooted tree with the root $v$. If for each $u\in V(T)\setminus\{v\}$, $d_T(u)\le2$, then
%\begin{equation}
%\psi(T)=\sum_{w\in C(v)}\psi(T_w)+\alpha(v),\notag
%\end{equation}
%where
\begin{center}
$\psi(T)=\left \{
  \begin{array}{ll}
   \sum\limits_{w\in C(v)}\psi(T_w)+1&,  \hbox{if $|C^2(v)|=0$ and $|C^1(v)|\le 1$ ;}\vspace{8pt}\\
   \sum\limits_{w\in C(v)}\psi(T_w)&, \hbox{otherwise.}
\end{array}
\right.$
\end{center}
\end{lem}
\pf Because $T_w$ is a path for each $w\in C(v)$, it is easy to determine $\psi(T_w)$ and $C^i(v)\cap C^j(v)=\emptyset$ for $i\neq j$. Note that $\sum_{w\in C(v)}\psi(T_w)\le\psi(T)\le\sum_{w\in C(v)}\psi(T_w)+1$. We consider the two cases.\par
\textbf{Case 1.} $|C^2(v)|=0$ and $|C^1(v)|\le 1$.\par
If $w\in C^0(v)$, then $T_w\cong P_n$ with $n\equiv 0(\bmod3)$. Let $F_w$ be a maximum dissociation set of $T_w$ such that $w\notin F_w$ ($F_w$ exists by Observation \ref{obs2.1}(b)). If $w\in C^1(v)$, then $T_w\cong P_n$ with $n\equiv 1(\bmod3)$ and let $F_w$ be a maximum dissociation set of $T_w$ such that $d_{T_w[F_w]}(w)=0$ ($F_w$ exists by Observation \ref{obs2.1}(c)). Now, let
\begin{equation}
F=\bigcup_{w\in C(v)}F_w \cup\{v\},
\end{equation}
then $F$ is a dissociation set of $T$ and $|F|=\sum_{w\in C(v)}\psi(T_w)+1$. Thus, $F$ is a maximum dissociation set of $T$ and $\psi(T)=\sum_{w\in C(v)}\psi(T_w)+1$.\par

\textbf{Case 2.} $|C^2(v)|\ge1$ or $|C^1(v)|\ge 2$.\par
Suppose, for a contradiction, that $\psi(T)=\sum_{w\in C(v)}\psi(T_w)+1$. Let $F$ be a maximum dissociation set of $T$. Then, $v\in F$ and $F\cap T_w$ is a maximum dissociation set of $T_w$ for each $w\in C(v)$. Let $F_w:=F\cap T_w$ for each $w\in C(v)$. If $w\in C^2(v)$, then $T_w\cong P_n$ with $n\equiv 2(\bmod3)$. By Observation \ref{obs2.1}(d), $w\in F_w$ and $d_{T_w[F_w]}(w)=1$. Thus, if $|C^2(v)|\ge1$, there is a 3-path in $T[F]$ that contains the vertex $v$, a contradiction.\par
If $w\in C^1(v)$, then $T_w\cong P_n$ with $n\equiv 1(\bmod3)$. By Observation \ref{obs2.1}(c), we have $w\in F_w$. Thus, if $|C^1(v)|\ge 2$, then there is a 3-path in $T[F]$ that contains the vertex $v$, a contradiction.

The proof is complete.\qed

\begin{theo}\label{the2.1}
Let $T$ be a rooted tree with the root $v$. If for each $u\in V(T)\setminus\{v\}$, $d_T(u)\le2$, then\\
(a) $v\in \mathcal{A}(T)$ if and only if $|C^2(v)|=0$ and $|C^1(v)|\le 1$;\\
(b) $v\in \mathcal{N}(T)$ if and only if $|C^2(v)|=2$ or $|C^1(v)|+|C^2(v)|\ge 3$.
\end{theo}
\pf (a) Necessity. Suppose, for a contradiction, that $|C^2(v)|\ge1$ or $|C^1(v)|\ge2$. Let $F=\bigcup_{w\in C(v)}F_w$, where $F_w$ is a maximum dissociation set of $T_w$ for each $w\in C(v)$. By Lemma \ref{lem2.1}, we have $\psi(T)=\sum_{w\in C(v)}\psi(T_w)$. Thus, $F$ is a maximum dissociation set of $T$ and $v\notin F$, which contradicts with $v \in \mathcal{A}(T)$.\par

Sufficiency. Suppose that $|C^2(v)|=0$ and $|C^1(v)|\le 1$. Then $\psi(T)=\sum_{w\in C(v)}\psi(T_w)+1$ by Lemma \ref{lem2.1} and the vertex $v$ is in all maximum dissociation sets of $T$. Thus, $v\in \mathcal{A}(T)$.

The proof of (a) is complete.

(b) Necessity.  Suppose, for a contradiction, that $|C^2(v)|\neq 2$ and $|C^1(v)|+|C^2(v)|\le 2$.

If $|C^2(v)|=0$ and $|C^1(v)|\le1$, then $v\in \mathcal{A}(T)$ by
(a), a contradiction.\par

If $|C^2(v)|=0$ and $|C^1(v)|=2$, then we assume $w_1,w_2\in
C^1(v)$. For each $w\in C^0(v)$, there exists a maximum dissociation
set $F_w$ of $T_w$ such that $w\notin F_w$ by Observation
\ref{obs2.1}(b). For $w_1\in C^1(v)$, we have $T_{w_1}-w_1  \cong
P_n$ with $n\equiv 0(\bmod3)$. Let $F_{w_1}$ be a maximum
dissociation set of $T_{w_1}-w_1$. Then $|F_{w_1}|=\psi(T_{w_1})-1$
by Observation \ref{obs2.1}(a). For $w_2\in C^1(v)$, there exists a
maximum dissociation set $F_{w_2}$ of $T_{w_2}$ such that
$d_{T_{w_2}[F_{w_2}]}(w_2)=0$ by Observation \ref{obs2.1}(c). Let
$F=\bigcup_{w\in C^0(v)}F_w \cup F_{w_1} \cup F_{w_2}\cup \{v\}$,
then $F$ is a dissociation set of $T$ and $|F|=\sum_{w\in
C(v)}\psi(T_w)$. By Lemma \ref{lem2.1}, $F$ is a maximum
dissociation set of $T$ and $v\in F$, a contradiction.\par

If $|C^2(v)|=1$ and $|C^1(v)|\le1$, then we assume $w_3\in C^2(v)$.
For $w_3\in C^2(v)$, we have $T_{w_3}-w_3 \cong P_n$ with $n\equiv
1(\bmod3)$. Let $F_{w_3}$ be a maximum dissociation set of
$T_{w_3}-w_3$, then $|F_{w_3}|=\psi(T_{w_3})-1$ by Observation
\ref{obs2.1}(a). For each $w\in C^0(v)$, let $F_w$ be a maximum
dissociation set of $T_w$ such that $w\notin F_w$. For $w\in
C^1(v)$, let $F_w$ be a maximum dissociation set of $T_{w}$ such
that $d_{T_{w}[F_{w}]}(w)=0$. Now let $F=\bigcup_{w\in
C(v)-\{w_3\}}F_w  \cup F_{w_3}\cup\{v\},$ then $F$ is a dissociation
set of $T$ and $|F|=\sum_{w\in C(v)}\psi(T_w)$. By Lemma
\ref{lem2.1}, $F$ is a maximum dissociation set of $T$ and $v\in F$,
a contradiction.\par

Sufficiency. Suppose for a contradiction that $F$ is a maximum dissociation set of $T$ and $v\in F$. For each $w\in C^0(v)$, $|F\cap T_w|\le \psi(T_w)$. If $w\in C^2(v)$, then $T_w\cong P_n$ with $n\equiv 2(\bmod3)$. Since $v\in F$, we have $|F\cap T_w|\le \psi(T_w)-1$. If $w\in C^1(v)$, then every maximum dissociation set $F_w$ of $T_w$ contains the vertex $w$. Since $v\in F$, there are at least $|C^1(v)|-1$ vertices $w$ in $C^1(v)$ such that $|F\cap T_w|\le \psi(T_w)-1$.
Hence, it is easy to check $|F|<\sum_{w\in C(v)}\psi(T_w)\leq\psi(T)$, a contradiction.

The proof of (b) is complete.\qed

%for any a maximum dissociation set $F_w$ of $T_w$, we have $w\in F_w$ and $d_{T_w[F_w]}(w)=1$. Since $v\in F$, for any $w\in C^2(v)$, $|F\cap T_w|\le \psi(T_w)-1$. If $w\in C^1(v)$, then $T_w\cong P_n$ with $n\equiv 1(\bmod3)$. By Observation \ref{obs2.1}(c), for any a maximum dissociation set $F_w$ of $T_w$, we have $w\in F_w$. Since $|F\cap C(v)|\le 1$,  there are at least $|C^1(v)|-1$ vertices in $C^1(v)$ such that $|F\cap T_w|\le \psi(T_w)-1$. Thus, if $|C^1(v)|=0$, then $|F|\le\sum_{w\in C(v)}\psi(T_w)-|C^2(v)|+1$; if $|C^1(v)|\ge1$, then $|F|\le\sum_{w\in C(v)}\psi(T_w)-|C^2(v)|-(|C^1(v)|-1)+1$. Since $|C^2(v)|=2$ or $|C^1(v)|+|C^2(v)|\ge 3$, it is easy to check $|F|<\sum_{w\in C(v)}\psi(T_w)<\psi(T)$, a contradiction.\qed

\section{Characterizations of $\mathcal{A}(T)$ and $\mathcal{N}(T)$}

A technique called \textit{pruning process} was introduced in \cite{Mynhardt1999}. Using the technique and Theorem \ref{the2.1}, we can characterize $\mathcal{A}(T)$ and $\mathcal{N}(T)$ for an arbitrary tree $T$.\par

Let $T$ be a rooted tree with the root $v$ and $u$ be a branch vertex of $T$ at maximum distance from $v$. It is easy to see that $|C(u)|\ge2$ and $d_T(x)\le2$ for each $x\in D(u)$. If $u\neq v$, we execute the following pruning process:\par
\begin{itemize}
\item if $|C^2(u)|\ge1$ or $|C^1(u)|\ge2$, then delete $D[u]$,
\item if $|C^2(u)|=0$ and $|C^1(u)|\le1$, then for all $w\in C(u)\setminus\{z\}$, delete $D[w]$, where $z$ is the vertex in $C^1(u)$ if $|C^1(u)|=1$ , otherwise $z$ is any one vertex in $C(u)$.
\end{itemize}
This step of pruning process is called a \textit{pruning of $T$ at $u$}.
%This step of pruning process in a tree $T_v$, where $D[u]$ or all but one child of $u$ together with their descendants are deleted to give a tree $T'$ such that $u\notin V(T')$ or $d_{T'}(u)=2$, is called a \textit{pruning of $T_v$ at $u$}.
%Let $T=(V, E)$ be a tree and $v\in V$ be an arbitrary vertex of $T$. Let $T_v$ be the rooted tree obtained from $T$ with the root $v$.
Repeat the above pruning process, finially we obtain a unique tree $\bar{T}$ called \textit{the pruning of $T$} such that $d_{\bar{T}}(u)\le2$ for each $u\in V(\bar{T})\setminus\{v\}$. We
will show that the root $v$ is in all maximum dissociation sets (or
in no maximum dissociation set) of $T$ if and only if it is in all
maximum dissociation sets (or in no maximum dissociation set) of the
pruning $\bar{T}$ of $T$.

\begin{lem}\label{lem3.1}
Let $T$ be a rooted tree with the root $v$ and $\bar{T}$ be the pruning of $T$. For every maximum dissociation set $\bar{F}$ of $\bar{T}$, there exists a maximum dissociation set $F$ of $T$ such that $v\in F$ if and only if $v\in\bar{F}$. Conversely, for every maximum dissociation set $F$ of $T$, there exists a maximum dissociation set $\bar{F}$ of $\bar{T}$ such that $v\in\bar{F}$ if and only if $v\in F$.
\end{lem}

\pf We prove the lemma by induction on $|B'(T)|$, where $B'(T)=\{u\in V(T)\setminus\{v\}:d(u)\ge3\}$. If $|B'(T)|=0$, then $T=\bar{T}$ and the result holds clearly. Suppose that when $|B'(T)|<k$, the lemma holds. Let $T$ be a tree with $|B'(T)|=k$ and $u$ be a vertex of $B'(T)$ at maximum distance from $v$. Let $T'$ be the tree obtained from $T$ by applying a pruning of $T$ at $u$. Thus, $\bar{T}$ is also the pruning of $T'$.

First, we show that for every maximum dissociation set $\bar{F}$ of $\bar{T}$, there exists a maximum dissociation set $F$ of $T$ such that $v\in F$ if and only if $v\in\bar{F}$. By the induction hypothesis, for every maximum dissociation set $\bar{F}$ of $\bar{T}$, there exists a maximum dissociation set $F'$ of $T'$ such that $v\in F'$ if and only if $v\in\bar{F}$.

We consider the following two cases.\par

\textbf{Case 1.} $|C^2(u)|=0$ and $|C^1(u)|\le 1$.\par
%Let $T'$ be a tree obtained from $T$ by deleting $D[w]$ for each $w\in C(u)-\{z\}$, where $z$ is the vertex in $C^1(v)$ if $|C^1(u)|=1$ , else $z$ is any one vertex in $C(u)$.

For each $w\in C(u)\setminus\{z\}$, $T_w\cong P_n$ with $n\equiv 0\ (\bmod3)$. Let $F_w$ be a maximum dissociation set of $T_w$ such that $w\notin F_w$ ($F_w$ exists by Observation \ref{obs2.1}(b)). Let $F=\bigcup_{w\in C(u)\setminus\{z\}}F_w \cup F'$. Clearly, $F$ is a maximum dissociation set of $T$. Since $v\in F'$ if and only if $v\in F$, we have $v\in F$ if and only if $v\in\bar{F}$.

\textbf{Case 2.} $|C^2(u)|\ge1$ or $|C^1(u)|\ge2$.\par

In this case, $T'=T-D[u]$. For each $w\in C(u)$, let $F_w$ be a maximum dissociation set of $T_w$. Let $F=\bigcup_{w\in C(u)}F_w \cup F'$. Since $|C^2(u)|\ge1$ or $|C^1(u)|\ge2$, by Lemma \ref{lem2.1}, $\psi(T_u)= \sum_{w\in C(u)}\psi(T_w)= \sum_{w\in C(u)}|F_w|$, which implies that $\bigcup_{w\in C(u)}F_w$ is a maximum dissociation set of $T_u$. Thus, $F$ is a maximum dissociation set of $T$. Since $v\in F'$ if and only if $v\in F$, we have $v\in F$ if and only if $v\in\bar{F}$.

Now, we show that for every maximum dissociation set $F$ of $T$ , there exists a maximum dissociation set $\bar{F}$ of $\bar{T}$ such that $v\in\bar{F}$  if and only if $v\in F$. Likewise, we consider the following two cases.\par

\textbf{Case 1.} $|C^2(u)|=0$ and $|C^1(u)|\le 1$.\par

For each $w\in C(u)\setminus\{z\}$, we have $T_w\cong P_n$ with $n\equiv 0(\bmod3)$. Let $F_w$ be a maximum dissociation set of $T_w$ such that $w\notin F_w$ ($F_w$ exists by Observation \ref{obs2.1}(b)).

For every maximum dissociation set $F$ of $T$, let $F'=F-\bigcup_{w\in C(u)\setminus\{z\}}D[w]$. Then $F'$ is a dissociation set of $T'$ and $v\in F$ if and only if $v\in F'$. We will prove that $F'$ is a maximum dissociation set of $T'$. Suppose, for a contradiction, that $F_1$ is a maximum dissociation set of $T'$ with $|F_1|>|F'|$. Let $F_2=\bigcup_{w\in C(u)\setminus\{z\}}F_w \cup F_1$. Then $F_2$ is a dissociation set of $T$ and
\begin{align}
|F_2|=|F_1|+\sum_{w\in C(u)\setminus\{z\}}|F_w| >|F'|+\sum_{w\in C(u)\setminus\{z\}}|F\cap D[w]| =|F|,
\end{align}
which contradicts to the fact that $F$ is a maximum dissociation set of $T$.
Thus $F'$ is a maximum dissociation set of $T'$. By the induction hypothesis, there exists a maximum dissociation set $\bar{F}$ of $\bar{T}$ such that $v\in F'$ if and only if $v\in\bar{F}$. Hence, there exists a maximum dissociation set $\bar{F}$ of $\bar{T}$ such that $v\in F$ if and only if $v\in\bar{F}$. \par

\textbf{Case 2.} $|C^2(u)|\ge1$ or $|C^1(u)|\ge2$.\par

For each $w\in C(u)$, let $F_w$ be a maximum dissociation set of $T_w$. By Lemma \ref{lem2.1}, $\psi(T_u)= \sum_{w\in C(u)}\psi(T_w)= \sum_{w\in C(u)}|F_w|$.

For every maximum dissociation set $F$ of $T$, let $F'=F-D[u]$. We will prove that $F'$ is a maximum dissociation set of $T'$. Suppose, for a contradiction, that $F_1$ is a maximum dissociation set of $T'$ with $|F_1|>|F'|$. Let $F_2=\bigcup_{w\in C(u)}F_w \cup F_1$. Then $F_2$ is a dissociation set of $T$ and \begin{align}
|F_2|=|F_1|+\sum_{w\in C(u)}|F_w| >|F'|+|F\cap D[u]| =|F|,
\end{align}
which contradicts with the fact that $F$ is a maximum dissociation set of $T$. Thus $F'$ is a maximum dissociation set of $T'$. By the induction hypothesis, there exists a maximum dissociation set $\bar{F}$ of $\bar{T}$ such that $v\in F'$ if and only if $v\in\bar{F}$. Thus, there exists a maximum dissociation set $\bar{F}$ of $\bar{T}$ such that $v\in F$ if and only if $v\in\bar{F}$. \par
We complete the proof.\qed

By Lemma \ref{lem3.1}, we can obtain the following corollary.
\begin{coro}\label{co4.1}
Let $T$ be a rooted tree with the root $v$ and $\bar{T}$ be the pruning of $T$, then $v\in\mathcal{A}(T)$ (or $\mathcal{N}(T)$) if and only if $v\in \mathcal{A}(\bar{T})$ (or $\mathcal{N}(\bar{T}))$.
\end{coro}

By Theorem \ref{the2.1} and Corollary \ref{co4.1}, the characterizations of $\mathcal{A}(T)$ and $\mathcal{A}(T)$ can be obtained immediately.

\begin{theo}
Let $T$ be a tree and $v$ be a vertex of $T$. Let $T_v$ be the rooted tree obtained from $T$ with the root $v$ and $\bar{T_v}$ be the pruning of $T_v$. Then

(a) $v\in\mathcal{A}(T)$ if and only if $|C^2_{\bar{T_v}}(v)|=0$ and $|C^1_{\bar{T_v}}(v)|\le1$;

(b) $v\in\mathcal{N}(T)$ if and only if $|C^2_{\bar{T_v}}(v)|=2$ or $|C^1_{\bar{T_v}}(v)|+|C^2_{\bar{T_v}}(v)|\ge3$.
\end{theo}

%\newpage
\section{A recognition algorithm}

In this section, we present a linear time recognition algorithm which can determine whether a
given vertex in a tree is in all (or no) maximum dissociation sets of the
tree. The recognition algorithm is described in detail as follows.
%
%By Theorem \ref{the2.1} and the pruning process, we can determine whether the root of a rooted tree is in all (or no)
%maximum dissociation sets. Further,

\begin{table}[htbp]
%\caption{\label{tab:test}11}
\begin{tabular}{lcl}
\toprule  \textbf{Recognition Algorithm} \\
\midrule
\textbf{Input}: a tree $T$ and a vertex $v\in T$. \\
\textbf{Output}: $v \in \mathcal{A}(T)$; or $v \in
\mathcal{F}(T)$; or $v \in
\mathcal{N}(T)$. \\

\quad 1. change the tree $T$ into a rooted tree by choosing the vertex $v$ as the root\\
\quad 2. compute the distance $d(v,u)$ from $v$ to each other vertex $u$\\
\quad 3. let $B = B(T) \setminus \{v\}$, where $B(T)$ is the set of branch vertices of $T$\\
\quad 4. \textbf{while} $B \neq \emptyset$ \textbf{do} \\
\quad    \quad 4.1. choose a vertex $u$ in $B$ such that $d(v,u)$ is maximum\\
\quad    \quad 4.2. \textbf{if} $|C^2(u)|\ge1$ or $|C^1(u)|\ge2$, \textbf{then}\\
\quad    \qquad \qquad $T\leftarrow T-D[u]$ and $B \leftarrow B \setminus \{u\}$\\
\quad    \quad 4.3. \textbf{else if} $|C^1(u)| = 1$, \textbf{then}\\
\quad    \qquad \qquad $T \leftarrow T-\bigcup\limits_{w \in C(u)\setminus\{z\}} D[w]$ and $B \leftarrow B \setminus \{u\}$, where $z$ is the vertex in $C^1(u)$\\
\quad    \quad 4.4. \textbf{else} \\
\quad    \qquad \qquad choose any one vertex $z$ in $C(u)$ and \\
\quad    \qquad\qquad $T\leftarrow T-\bigcup\limits_{w \in C(u)\setminus\{z\}} D[w]$ and $B \leftarrow B \setminus \{u\}$\\
\quad 5. \textbf{If} $|C^2(v)|=0$ and $|C^1(v)|\le 1$, \textbf{then}\\
\quad    \qquad output $v\in \mathcal{A}(T)$\\
\quad 6. \textbf{else if} $|C^2(v)|=2$ or $|C^1(v)|+|C^2(v)|\ge 3$, \textbf{then}\\
\quad    \qquad output $v\in \mathcal{N}(T)$\\
\quad 7. \textbf{else}\\
\quad    \qquad  output $v \in \mathcal{F}(T)$\\
\bottomrule
\end{tabular}
\end{table}

For a rooted tree, we can use the breadth-first search algorithm to
find the distance from the root to each other vertex in linear time.
Step 4 is the pruning process of the rooted tree $T$ and can be
executed in linear time. Thus, the runtime of the recognition
algorithm is linear. Thus, using the recognition algorithm we can
find all vertices contained in all (or no) maximum dissociation sets
of a tree of order $n$ in $O(n^2)$ time.

\section*{Acknowledgments}
\noindent The work is supported by Research Foundation for Advanced Talents of Beijing Technology and Business University. 

\bibliographystyle{unsrt}
%\bibliography{acl}

\end{document}